\documentclass[aps,apl,reprint,groupedaddress, showkeys,showpacs]{revtex4-1}
\usepackage[utf8]{inputenc}
\usepackage{amsmath,amssymb}
\usepackage{graphicx,float}
\makeatletter
\usepackage{dsfont}
\usepackage{multirow}
\usepackage{url}
\makeatother

\begin{document}

\title{General expression for the component size distribution in infinite configuration networks}
\author{Ivan Kryven} 
\email{i.kryven@uva.nl}
\affiliation{University of Amsterdam, PO box  94214, 1090 GE, Amsterdam, The Netherlands} 

\begin{abstract}
In the infinite configuration network the links between nodes are assigned randomly with the only restriction that the degree distribution has to match a predefined function. This work presents a simple equation that gives for an arbitrary degree distribution  the corresponding size distribution of connected components. This equation is suitable for fast and stable numerical computations up to the machine precision. The analytical analysis reveals that the asymptote of the component size distribution is completely defined by only a few parameters of the degree distribution: the first three moments, scale and exponent (if applicable). When the degree distribution features a heavy tail, multiple asymptotic modes are observed in the component size distribution that, in turn, may or may not feature a heavy tail.
\end{abstract}



\maketitle

\section{Introduction}
Random graphs provide models for complex networks, and in many cases, real-world networks has been accurately described by such models\cite{barabasi1999,newman2002,vazquez2003,kryven2016b}. 
Within the scope of random graph models one finds: Erd\H{o}s-R\'{e}nyi model, Barab\'{a}si-Albert model\cite{barabasi1999}, node copying model\cite{bhat2016}, small world network\cite{newman1999}, configuration network\cite{molloy1998a,newman2010book} and many others.
In the configuration network $N$ nodes are assigned pre-defined degrees.
The edges connecting these nodes are then considered to be random, and every distinct configuration of edges that satisfies the given degree sequence is treated as a new instance of the network in the sense of random graphs. Interesting properties emerge when the number of nodes, $N,$ approaches infinity, or at  the so-called thermodynamic limit\cite{chung2003,Kryven2016a}. In this case, the infinite degree sequence, which provides the only input information for the model, is equivalent to the frequency distribution of degrees, $u(k), \;k=1,2,\dots,$ i.e. the probability that a randomly chosen node has degree $k$.

Component-size distribution, $w(n),$ denotes probability that a randomly chosen node is part of a connected component of finite size $n.$
 Connected components in the infinite configuration network can be of finite or infinite size.  Molloy and Reed\cite{molloy1995} showed  that if an infinite component exists then it is the only infinite component with probability 1. Hence, the infinite component is referred to as the giant component.

Depending upon a specific context behind the network, the component size distribution may summarise an important feature of the modelled system.
In polymer chemistry, for example, the infinite configuration network is used as a toy model for hyper-branched and cross-linked polymers. 
In this context, the component size distribution predicts viscoelastic properties of the material while the emergence of the giant component is interpreted as a phase transition from liquid to solid state of the soft matter\cite{kryven2016b,kryven2014c}.
Since connected components are closely related to clusters in bond percolation processes, 
the distribution of component sizes can be used to model outbreaks for SIR epidemiological processes\cite{newman2002}.
In linguistics, component size distribution of the sentence similarity graph is an important tool when studying structure of natural languages  \cite{biemann2011}. This brief list of application cases is far from being exhaustive.
Despite the vast applications, the  empirical component size distribution is hard to measure precisely unless the whole topology of the network is known.
On another hand, empirical observations on the degree distribution, $u(k),$ are much easier to perform.

Ref.\cite{molloy1995} provides an elegant criterion that connects moments of $u(k)$ to the fact that the network contains the giant component. 
A somewhat deeper question further in this direction reads: providing $u(k)$ is given, what is the component size distribution, $w(n)?$ 
In Ref.\cite{Newman2001} Newman et al. showed that the component size distribution can be recovered by a numerical algorithm that involves solving a fixed point problem followed by an inversion of a generation function. 
Such algorithm demonstrates that indeed $u(k)$ and $w(n)$ can be put into a correspondence, however, it becomes computationally infeasible for large values of $n$. This numerical issues aries due to ill-posedness of the numerical generating function inversion.

On another hand, the component size distribution has been analytically resolved only for a limited number of partial cases of $u(k)$\cite{Newman2007}. Within the scope of analytically solvable cases, only the Yule-Simon degree distribution features a heavy tail, that is to say it decays proportionally to an algebraic function, $n^{-\beta},\;\beta>0$ at  large $n$.
At the same time, the heavy-tailed  (or scale-free) distributions are commonly observed in the empirical data collected from many real-world networks\cite{yook2002,ravasz2003,eguiluz2005,fu2008}. Empirically observed exponents vary in a broad range. Some studies report degree exponents that are as small as $\beta=0.81$ in the case of the Internet topology\cite{adamic2002} and $\beta=1$ in social networks \cite{timar2016}. On the opposite side of this spectrum, one finds exponent $\beta=5$ in the generalisation of preferential attachment model \cite{dorogovtsev2000}.

The only asymptotic analysis available for component size distribution in the configuration network states that for large $n,$ $w(n)$ is either proportional to $n^{-3/2}$ or to $e^{-C n},$ where $C>0$ is a constant \cite{Newman2001}. The current paper uncovers new asymptotic modes for $w(n)$ that emerge only when the degree distribution features a heavy tail. The paper shows that for an arbitrary $u(k),$ $w(n)$ can be expressed as a finite sum. 
In practice, this sum can be stably computed up to the machine precision in the cost of $O(n^2)$ multiplicative operations. Finally, the paper discusses how a finite cutoff introduced in the degree distribution reflects on the distribution of component sizes.

\section{Component size distribution by Lagrange inversion }
It has been noticed that all components in the infinite configuration model are locally tree-like.  Using this fact as a departure point, Newman et al.\cite{Newman2001} showed that the degree distribution can be put into a correspondence to the component size distribution by applying the generation-function (GF) formalism. 
Here, by a GF of $u(k),\;   \sum_{k=0}^\infty u(k)=1,$ we refer to the series,
\begin{equation}
\label{eq:transform}
U(x) =  \sum_{k=0}^\infty u(k) x^k, \; x\in\mathbb{C}, \;|x| \leq 1.
\end{equation}
According to the approach presented in Ref.~\cite{Newman2001}, the generating function for the component size distribution, $W(x)$ is found as a solution of the following system of functional equations,
\begin{align}
\label{eq:Wa}
W(x) = x U[W_1(x)],\\
\tag{\ref{eq:Wa}$'$}
\label{eq:Wb}
W_1(x) = x U_1[W_1(x)].
\end{align}
where $U(x)$ is the GF of $u(k)$, and $U_1(x)$ is the GF for the excess degree distribution 
\begin{equation}\label{eq:defu1}
u_1(k) = \frac{k+1}{\mu_1} u(k+1),
\end{equation}
where $\mu_1= \sum_{k=1}^\infty k u(k).$
Similarly to combinatorial tree-counting problems, Eq.~\eqref{eq:Wa} can be solved  
by applying the Lagrange inversion formula\cite{bergeron1998}. 
The original formulation of the Lagrange inversion principle is as follows. Suppose, $A(X),R(x),$ are such formal power series that $A(x)=x R[A(x)]$ then
for an arbitrary formal power series $F(x),$ the coefficient of power series $F[A(x)]$ at $x^n$ reads as,
\begin{equation}
\label{eq:LagrangeInversion_1d}
[x^n]F[A(x)] = \frac{1}{n}[t^{n-1}]F'(t)R^n(t),\; n>0.
\end{equation}
Here $[t^{n-1}],$ as being the inverse operation to the GF transform \eqref{eq:transform}, refers to the coefficient at $t^{n-1}$ of the corresponding power series. 
By substituting $A(x)=W_1(x)$, $F(x)= U(x)$ and applying Eq.~\eqref{eq:Wa} one transforms the left hand side of \eqref{eq:LagrangeInversion_1d},
$
[x^n]F[A(x)]=[x^n]U[W_1(x)]=[x^{n-1}]W(x)=w(n-1).
$
 Further on, the right hand side of \eqref{eq:LagrangeInversion_1d} is transformed  by substituting $R(x)=U_1(x)$ and realising that, according to the definition \eqref{eq:defu1}, $U^{'}(x)=\mu_1 U_1(x).$ Now, we are ready to write an expression for $w(n),$  even though we have no explicit expression for generating function $W(x)$ itself, 
\begin{multline}
\label{eq:WviaLngrange}
w(n)=\frac{1}{n-1}[t^{n-2}]U^{'}(x) U_1(x)^{n-1}=\\
\frac{\mu_1}{n-1}[t^{n-2}] U_1(x)^{n},\; n>1.
\end{multline}
A similar equation was also derived in Ref.\cite{Newman2007} by means of different reasoning.
In principle, Eq.~\eqref{eq:WviaLngrange} provides enough information to analytically recover the component size distribution for a few special cases of the degree distribution\cite{Newman2007}. 
In practice, however, the main difficulty when applying Eq.~\eqref{eq:WviaLngrange} is that the equation employs the inverse GF transform, $[t^{n}]$, which limits the choices one has when searching for an exact solution or performing numerical computations. 
With this in mind, one may rewrite \eqref{eq:WviaLngrange} so that the new expression does not involve the GF concept at all. It turns out that the only reason why Eq. \eqref{eq:WviaLngrange} utilises the GF formalism is that it provides means for \emph{convolution power}.

The convolution of two distributions,  $f(k)*g(k), \;k>0$ is defined as a binary multiplicative operation, $$f(k)*g(k) =\sum\limits_{i+j=k}f(i)g(j),$$ where the summation is performed over all non-negative ordered couples $i,j$ that sum up to $k$. This sum contains exactly $k+1$ of such couples.
 In this paper, the order of operations is chosen in such a way that the point-wise multiplication precedes convolution, for instance, $f(k)*k g(k)=  f(k)*[k g(k)]$. The convolution can be inductively extended to the $n$-fold convolution, or  \emph{ the convolution power}, 
 \begin{equation}\label{eq:conpower}
 f(k)^{*n}=f(k)^{*n-1}*f(k),
 \end{equation}
  where $f(k)^{*0}\equiv 1 $ by the definition.
It can be shown that the convolution power can be expanded into a sum of products,  
\begin{equation}\label{eq:sum}
f(k)^{*n}=\sum\limits_{\stackrel{k_1+\dots +k_n=k}{k_i \geq 0}}\prod\limits_{i=1}^n f(k_i).
\end{equation}
 The convolution has a peculiar property in respect to the GF transform. If $F(x),G(x),U(x)$ are GFs for $f(k),g(k),$ and $u(k)=f(k)*g(k)$ then $U(x)=F(x)G(x)$. Furthermore, if $U(x)$ is GF for $u(k)$ then $U(x)^n$ generates $u(k)^{*n}.$
By exploiting this relation one immediately reduces Eq.~\eqref{eq:WviaLngrange} to, 
\begin{equation}
\label{eq:Lagrange1d}
w(n)=\begin{cases}
\frac{\mu_1}{n-1} u_1^{*n}(n-2),& n>1, \\
u(0) & n=1.
\end{cases}
\end{equation}
Here, the value of $w(0)$ is derived directly from the formulation of the problem: nodes with degree zero are also components of size one.
This simple equation is ready to be used: by combining \eqref{eq:Lagrange1d} and the definitions \eqref{eq:defu1},\eqref{eq:sum} one may directly expresses the values of the component size distribution in terms of $u(k):$ 
for $n>1,$
\begin{equation}
\label{eq:Lagrange1d*}
w(n)=
\frac{[ku(k)]^{*n}(2n-2)}{(n-1)\mu_1^{n-1}},\;n>1 .
\end{equation}
 For example,  first five values of $w(n)$ read as,
\begin{equation*}
\begin{aligned}
w(1) = & u(0),\\
w(2) = &\frac{1}{\mu_1} u(1)^2,\\
w(3) = &\frac{3}{\mu_1^2} u(1)^2u(2),\\
w(4) = &\frac{4}{\mu_1^3} u(1)^2[2 u(2)^2 + u(1) u(3)],\\
w(5) = &\frac{5}{\mu_1^4} u(1)^2[4u(2)^3 + 6 u(1) u(2) u(3) + u(1)^2 u(4)].
\end{aligned}
\end{equation*}
The number of terms in this expansion increases rapidly with $n$. That said, the formula \eqref{eq:Lagrange1d} can be easily readjusted for numerical computations. Namely, one can use Fast Fourier Transform (FFT) to compute the convolution powers, $ u_1^{*n}(k)=\mathcal{F}^{-1}[\mathcal{F}[u_1(k)]^n].$ In this case, $O(N^2)$ multiplicative operations is sufficient to compute all values of $w(n), n\leq N$. Alternatively, if $w(n-1)$ is known $w(n)$ can be found in the cost of $O(n \log n ).$ Besides FFT, there are algorithms that are specifically designed for fast approximation of convolution powers, such as projection onto basis functions that are invariant under convolution\cite{kryven2013a}.

\begin{figure}[htbp]
\begin{center}
\includegraphics[width=0.4\textwidth]{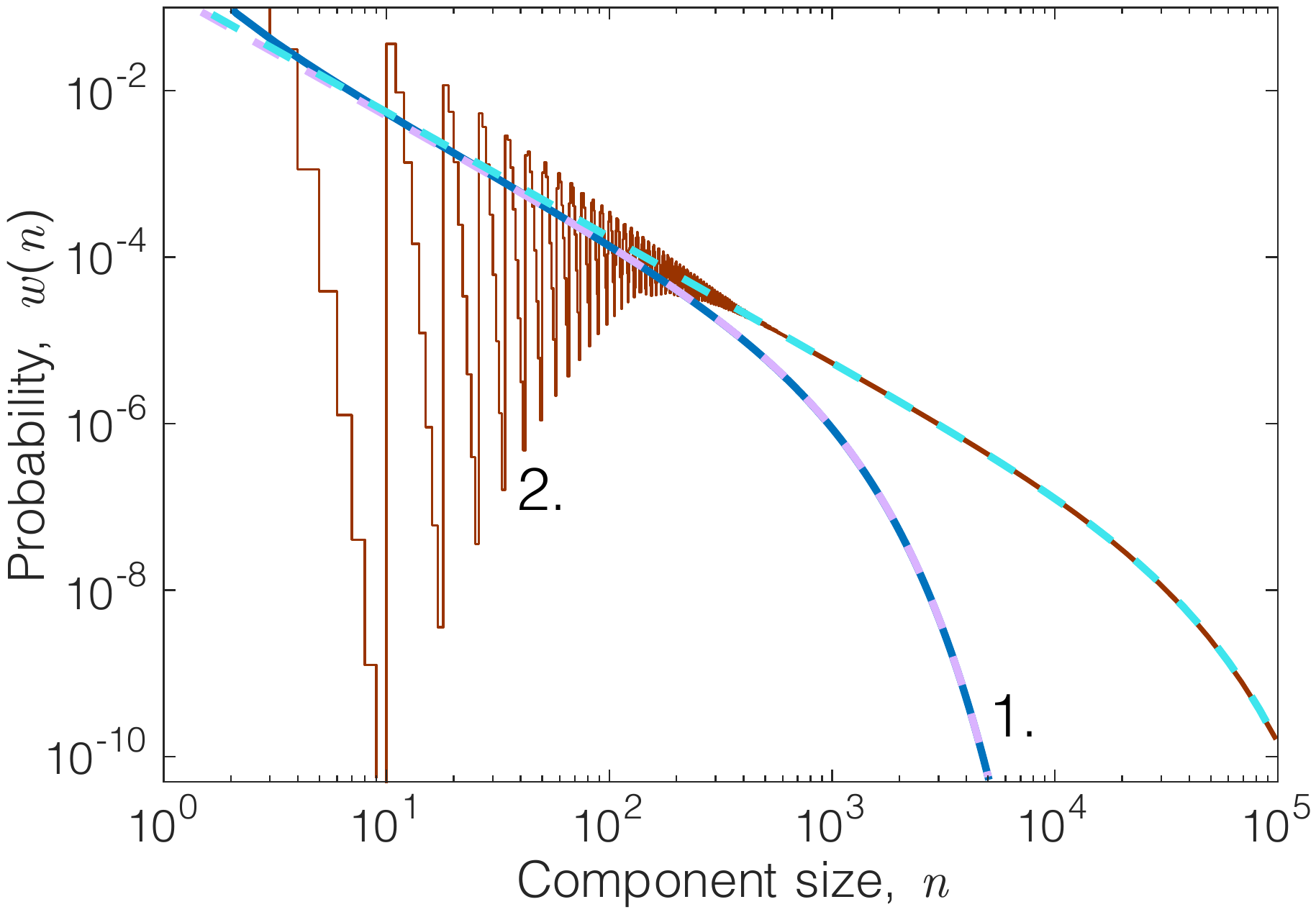}
\caption{  Examples of component size distributions (solid lines) that feature fast (\emph{1}) and slow (\emph{2}) convergence to their asymptotes (dashed lines). Both asymptotes are covered by Case A, Table~\ref{tab:asymptotes}. \emph{1})  $u(k)=C e^{-1.05k},$ all three: the analytical expression  (see Table~\ref{tab:ex1}), numerical values (according to Eq.~\eqref{eq:Lagrange1d}) and the asymptote practically coincide. \emph{2}) $u(k)$ is non-zero in three points $u(1)=0.97,\; u(2)=0.015,\;u(10)=0.015;$ the component size distribution features oscillations before it converges to the asymptote. }
\label{fig:zigzag}
\end{center}
\end{figure}

Analytic formulas for convolution powers (sometimes also referred to as \emph{compositas} \cite{kruchinin2015}), were covered by literature for many elementary functions\cite{ma2004,ma2002}. Convolution powers of $u_1(k)$ can also be found analytically by applying discrete functional transforms, for instance Z-transform and discrete Fourier transform. A few examples of such results are given in Table \ref{tab:ex1}. Focusing on one of them, the first curve in Figure~\ref{fig:zigzag} demonstrates that both analytical and numerical results for the exponential degree distribution coincide.
\begin{table}[htp]

\begin{center}
\begin{tabular}{ll}
Degree Distribution, $\;\;\;\;\;\;\;\;$ & Component Size Distribution\\
\hline
\multicolumn{2}{l}{Exponential distribution}\\
$C e^{-\lambda k}$
&
$\frac{(1-e^{-\lambda})^{2n-1}}{e^{\lambda(n-1)}}\frac{\Gamma (3 n-2)}{\Gamma (n) \Gamma (2 n)}$\\

\multicolumn{2}{l}{Geometric distribution}\\
$  (1 - p)^{k - 1} p$ &
$ (1-p)^{n-2} p^{2 n-1} \frac{ \Gamma (3 n-2) }{\Gamma (n) \Gamma (2 n)}$\\
\multicolumn{2}{l}{Binomial distribution}\\
\multirow{2}{*}{$ \binom{ k_{ \text{max}} }{k} (1 - c)^{k_\text{max}-k}  c^k $}&
$\frac{1}{n-1} \binom{ n k_{ \text{max}}-n}{n-2}\times  $ \\
 &$ \;\;\;\;\;\;\; (1 - c)^{ n k_{\text{max}}  - 2n-2}   c^{n-2}  $
\end{tabular}
\end{center}
\caption{Exact expressions for component size distributions in configuration network as evaluated with Eq.~\eqref{eq:Lagrange1d} via Z-Transform.}
\label{tab:ex1}
\end{table}%

\section{Asymptotic analysis }
The format of Eq.~\eqref{eq:Lagrange1d}, naturally suggests a straightforward way to perform an asymptotic analysis for $n \to \infty.$ One may view $u(k)$ as a probability mass function PMF (or alternatively discrete probability density function) of some discrete random variables $k_i.$  Recall the following property of convolution powers: if i.i.d. random variables $k_i$  have PMF $u_1(k)$ then $u_1^{*n}(k)$ gives the PMF for the sum  $k_1+k_2+\dots+k_n$. The central limit theorem (CLT)  gives an estimate for this sum as $n\to\infty,$
and the idea is now to obtain the asymptotes of $w(n)$ by applying CLT to the definition \eqref{eq:Lagrange1d}.
\subsection{Light-tailed degree distributions}
First, let us assume that  distribution  $u(k)$ decays faster than algebraically, that is 
\begin{equation}\label{eq:conds}
u(k) =o(k^{-\beta}),\; \beta>2,\; k \to\infty,
\end{equation}
which is also equivalent to $u_1(k) =o(k^{-\beta+1})$.
 Then according to CLT,  $u_1^{*n}(k)$ approaches the normal distribution, $u_1^{*n}(k) \stackrel{d}{\to} (\sqrt{n} \sigma)^{-1}\mathcal{N}(\frac{k-n M}{\sqrt{n} \sigma},0,1 ),$  when $n \to \infty,$ where  $M = \sum_{k=1}^{\infty} k u_1(k)$ and $\sigma^2 = \sum_{k=1}^{\infty}k(k-M)^2 u_1(k)<\infty$ denote the mean value and variance of $u_1(k)$.
The normal distribution can now replace $u_1^{*n}(k)$ in \eqref{eq:Lagrange1d}, which yields the asymptote for the component size distribution,
\begin{equation}
\label{eq:w_asm1}
w(n) \sim \frac{\mu_1 e^{  -\frac{  ( n (1 - M) -2 ) ^ 2 }{ 2 n  \sigma^2 }  }}{(n-1) \sqrt{2 \pi n \sigma^2}   } , \text{ as } n \to \infty.   
\end{equation}
Quantities $M,\sigma^2$ are directly expressible in terms of moments of degree distribution $u(k),$
\begin{equation}
\begin{aligned}\label{eq:M_sigma}
M   =&\sum_{k=1}^{\infty} k  u_1(k)= \frac{1}{\mu_1}\sum_{k=1}^{\infty} (k^2-k) u(k)=\frac{\mu_2-\mu_1}{\mu_1},\\
\sigma^2  =&\sum_{k=1}^{\infty}k(k-M)^2 u_1(k)= \\& \frac{1}{\mu_1}\sum_{k=0}^{\infty} k (k- M-1 )^2  u(k)= 
\frac{\mu_3 \mu_1-\mu_2^2}{\mu_1^2},
\end{aligned}
\end{equation}
where
$$\mu_i=\mathbb{E}[k^i]=\sum_{k=1}^{\infty}k^i u(k), i=1,2,\dots$$
Finally, substituting the expressions \eqref{eq:M_sigma} into \eqref{eq:w_asm1}  and gives the final version of the asymptote,
 \begin{equation}
\label{eq:w_asm2}
w(n) \sim \frac{\mu_1^2  
n^{-3/2}e^{-\frac{(\mu_2  -2 \mu_1 )^2}{2  (\mu_1 \mu_3 -\mu_2^2  )}n}
   }{ \sqrt{2 \pi(\mu_1 \mu_3-\mu_2^2)}}, \text{ as } n \to \infty. 
 \end{equation}
Two examples of component size distributions that converge with various rates to their asymptotes are given in Figure~\ref{fig:zigzag}.
Peculiarly, the only information on $u(k)$ that is contained in the asymptote definition \eqref{eq:w_asm2} is the first three moments $\mu_1,\mu_2,\mu_3$.
 Furthermore, depending upon the value of $\theta =\mu_2-2 \mu_1,$ the asymptotic expression \eqref{eq:w_asm2} switches between the \emph{two modes}: it either decays exponentially as $O(e^{-A n}),$ when $\theta \neq 0$, or it decays as an algebraic function, $O(n^{-3/2}),$ when $\theta=0$ (see also Table~\ref{tab:asymptotes}, Case A). The last equality is the well-known giant component criterion,
\begin{equation}
\label{eq:MolloyCriterion}
\mu_2-2 \mu_1=0.
\end{equation}
The criterion \eqref{eq:MolloyCriterion} was obtained by Molloy and Reed\cite{molloy1995} by means of a different reasoning. In Ref \cite{molloy1995}, the authors prove that $\theta>0$ implies existence of the giant component in the configuration network, whereas $\theta<0$ implies non-existence of this component. In Ref. \cite{Newman2001}, it was hypothesised that the $-3/2$ exponent is universal and must hold for all degree distributions at the critical point $\theta=0$. We will see now that when the condition \eqref{eq:conds} fails to hold,  distinct from $-3/2$ exponents may also appear in the asymptotic of $w(n)$.

\subsection{Heavy-tailed degree distributions}
\begin{figure}[htbp]
\begin{center}
\includegraphics[width=0.4\textwidth]{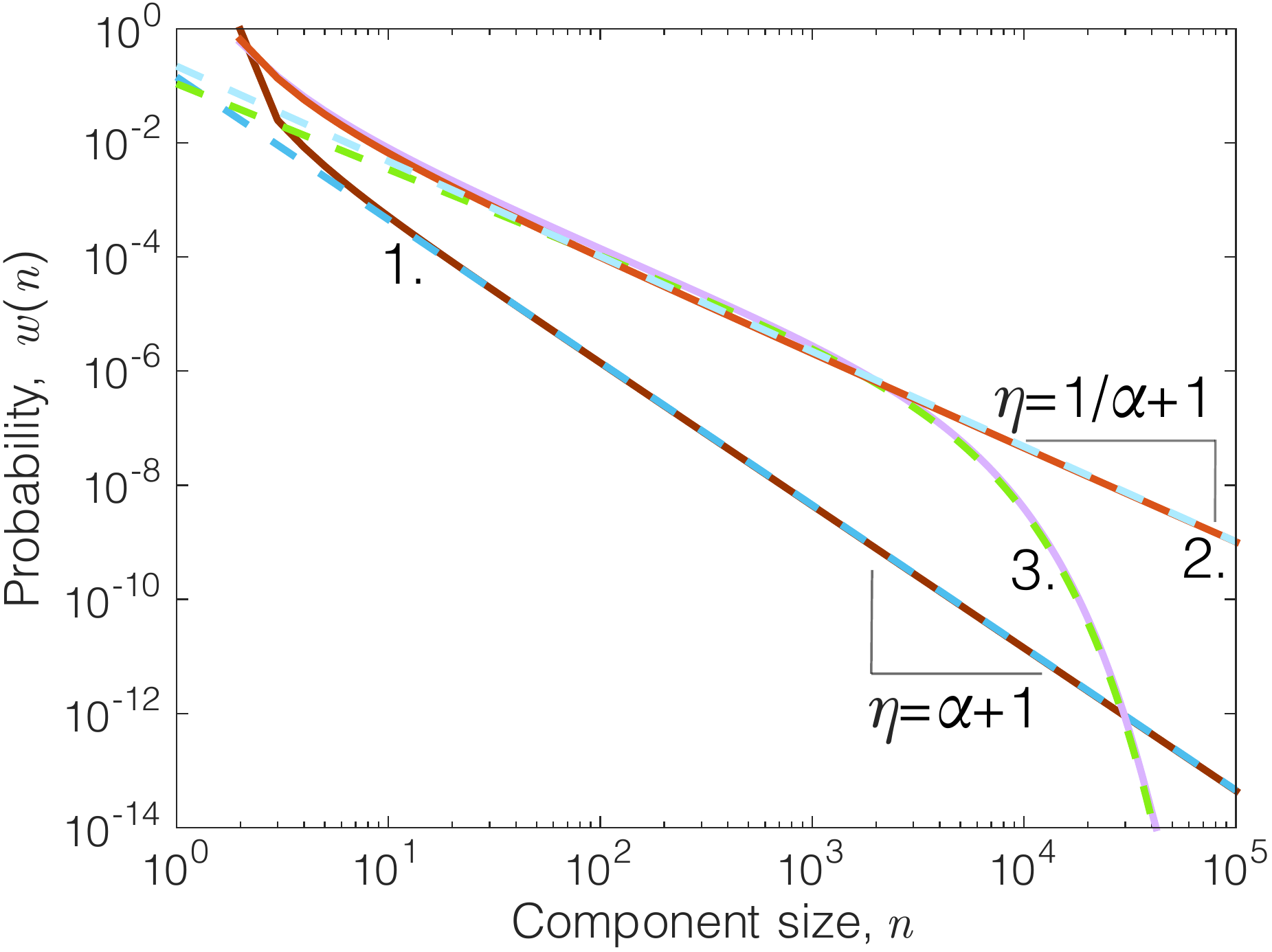}
\caption{ Component size distributions (\emph{solid lines}) and their asymptotes (\emph{dashed lines}) as obtained for degree distributions with exponent $\beta=3.5$ ($\alpha=1.5$) and various values of scale parameter. Three distinct asymptotic modes are illustrated:
 \emph{1}) $s = 0.066:$ $\theta<0,$        $\eta=\frac{5}{2}$; 
 \emph{2}) $s = 0.644:$ $\theta\approx0$,  $\eta=\frac{5}{3}$; 
 \emph{3}) $s = 0.8:$   $\theta>0$; 
   }
\label{fig:a15}
\end{center}
\end{figure}

Suppose that, on the contrary to  the condition \eqref{eq:conds},  degree distribution $u(k)$ features a heavy tail,  
\begin{equation}
\label{eq:tail_cond}
u(k)\sim s k^{-\beta},\; \beta>2,\;k\to \infty,
\end{equation}
which is equivalent to $u_1(k)\sim s k^{-\alpha-1},\; \alpha={\beta-2}>0$, $k\to \infty.$
It turns out that exponent $\alpha$ and  the scale $s,$ together with the moments $\mu_1,\mu_2,\mu_3$ provide enough information to generalise the asymptote \eqref{eq:w_asm2} for  the case of heavy-tailed degree distributions.
Suppose  $0<\alpha \leq 2.$
In terms of $u(k)$ moments this condition casts out as $\mu_3 = \infty.$ As follows from Gnedenko and Kolmogorov's generalisation of CLT\cite{gnedenko1954} the mass density distribution for $u_1^{*n}(k)$ approaches the stable law, 
\begin{equation}\label{eq:law}
u_1^{*n}(k) \stackrel{d}{\to} \frac{1}{\gamma(n)} G^A\left(\frac{k-\mu(n)}{\gamma(n)},\alpha,1\right),\;n\to\infty.
\end{equation}
 Here, we use the notation of Uchaikin \& Zolotarev \cite{uchaikin1999} which includes: \emph{exponent parameter} $\alpha,$ the \emph{location parameter}
\begin{equation}
\label{eq:Zol_mu}
\mu(n)=\begin{cases}
n \frac{\mu_2-\mu_1}{\mu_1}, & \alpha  >1,\\
s n \ln n, & \alpha  =1,\\
0, &0<\alpha<1,
\end{cases}
\end{equation}
and the \emph{scale parameter}
\begin{equation}
\label{eq:gamma}
\gamma(n)=\begin{cases}
\sqrt{s n \ln n}, & \alpha  =2,\\
\sqrt{\pi s}[2 \Gamma(\alpha)\sin \frac{\alpha \pi}{2}]^{-1/\alpha}n^{1/\alpha}, & \alpha \in (0,1)\cup (1,2)\\
\frac{\pi n s}{2}, & \alpha=1.
\end{cases}
\end{equation}
No general analytical expression  is known for $G^A(x,\alpha,1),$ and the stable law is defined via its Fourier transform
\begin{equation}
\mathcal{F}[{G^A(x,\alpha,1)}] = 
\begin{cases}
e^{-x^\alpha - i\tan \frac{\pi \alpha }{2}}, & \alpha \in (0,1)\cup (1,2],\\
e^{-x^\alpha + i\frac{2 \alpha }{\pi}}, & \alpha =1.
\end{cases}
\end{equation}
Consider the case when $1<\alpha<2.$ According to \eqref{eq:law} the point in which the stable law is evaluated, $x(n)=\frac{n-\mu(n)}{\gamma(n)}$, approaches  positive or negative infinities depending upon the sign of $\theta=\mu_2-2\mu_1.$ Indeed, as $n\to \infty,$
\begin{equation}
\label{eq:point_asympt}
x(n) \to \begin{cases}
+\infty ,& \theta<0,\\
0,& \theta=0,\\
-\infty ,& \theta>0.
\end{cases}
\end{equation}
For these values of $\alpha,$ function $G^A(x,\alpha,1)$ is non-zero on $(-\infty,+\infty).$ If $x(n)\to \infty,$ the function features an algebraic decay, whereas if $x(n)\to -\infty$ the decay is exponential. Therefore, the limiting value switching that takes place in   \eqref{eq:point_asympt} may reflect on the asymptotic behaviour of $u_1^{*n}(n).$ To give a precise answer one has to consider series expansions of $G^A(x,\alpha,1)$ around the points of interest, $x\in\{-\infty,0,+\infty\}$. We use here the leading terms of these series\cite{uchaikin1999}, 
\begin{multline}
\label{eq:lead_term}
G^A( x, \alpha,  1 ) =\\
 \begin{cases} 
\frac{\Gamma \left(1+\frac{1}{\alpha }\right) \sin \frac{\pi }{\alpha }}{\pi }+ O(x), & x\to 0\\
\frac{\Gamma (\alpha +1) x^{-\alpha -1}}{\Gamma (2-\alpha ) \Gamma (\alpha -1)}+O(x^{-2\alpha-1}),& x\to \infty,\\
\frac{e^{-(\alpha -1) \left(\frac{x}{\alpha }\right)^{\frac{\alpha }{\alpha-1}}} \left(\frac{x}{\alpha }\right)^{\frac{1}{2} \left(\frac{1}{\alpha -1}-1\right)}}{\sqrt{2 \pi  \alpha(\alpha -1) }} 
[1+ O(x^{-\frac{\alpha-1}{\alpha}}) ],& x\to -\infty,
\end{cases}
\end{multline}     
By replacing  the expression for the limiting distribution \eqref{eq:law} with the leading terms given in \eqref{eq:lead_term} one obtains the asymptotes for \eqref{eq:Lagrange1d}. This time, the asymptote has \emph{three modes}: depending upon the value of $\theta,$ it either features a heavy tail with exponent $-\alpha-1,$ a heavy tail with exponent $-\frac{1}{\alpha}-1,$ or an exponential decay, as shown in Table~\ref{tab:asymptotes}, Case D. 
A few examples of such asymptotic modes for a heavy-tailed degree distribution 
\begin{equation}\label{eq:degree_example}
u(k)=
\begin{cases}
C & k=1,\\
s(\beta-2) k^{-\beta}&k>1
\end{cases}
\end{equation}
are given in Fig.~\ref{fig:a15}.
The degree distribution \eqref{eq:degree_example} is defined by two parameters: exponent $\beta$ and scale $s$; whereas the constant $C$ is such that the total probability is normalised, $\sum_k u(k) =1$.

 When $\alpha=2$, the behaviour of $\frac{n-\mu(n)}{\gamma(n)},$  $n\to \infty$ is identical to \eqref{eq:point_asympt}, but the expression for $\gamma(n)$ is different and the series expansions \eqref{eq:lead_term} lead to somewhat different asymptotes, see Table~\ref{tab:asymptotes}, Case C. 
 
\begin{figure}[h]
\begin{center}
\includegraphics[width=0.4\textwidth]{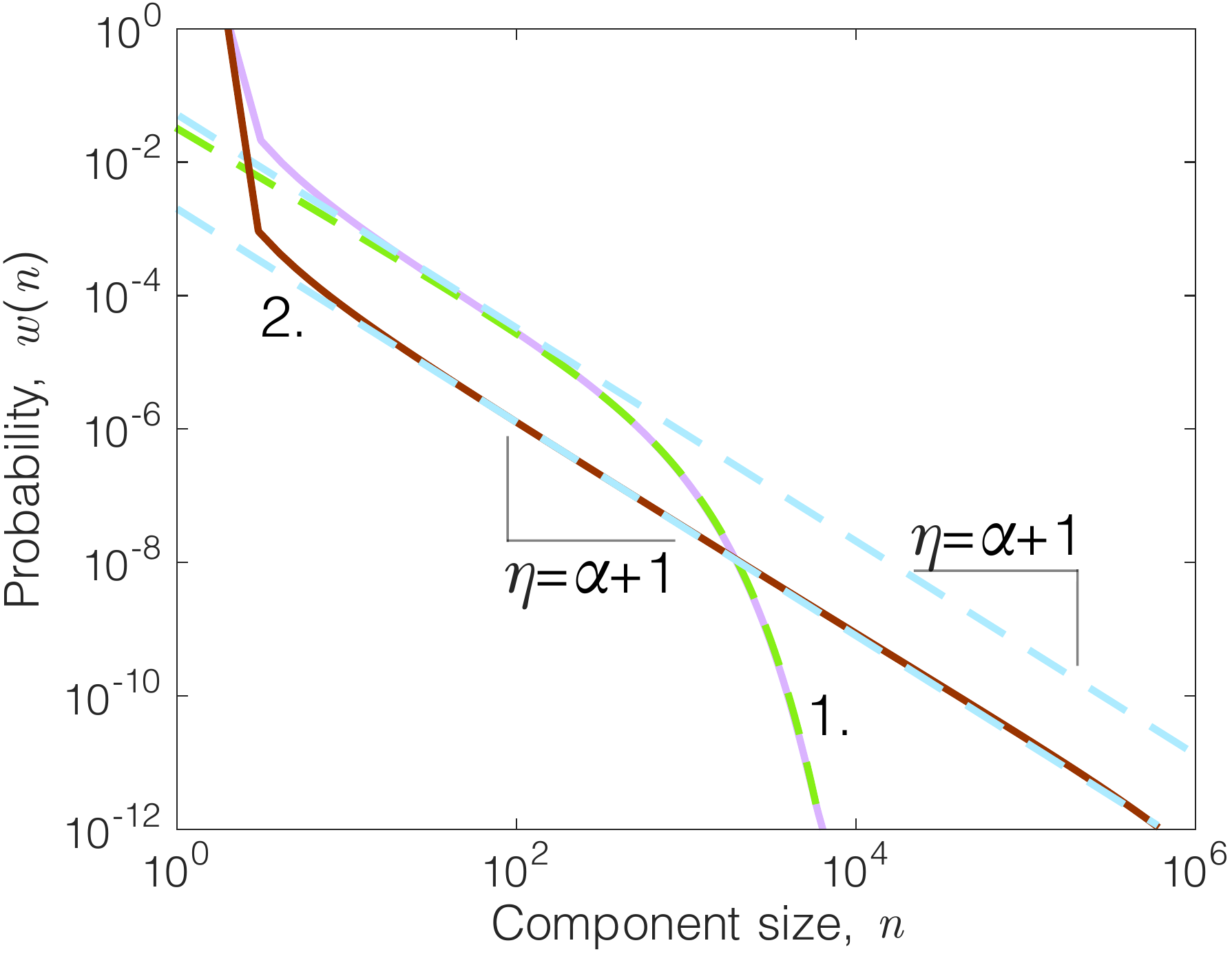}
\caption{ 
Component-size distributions (\emph{solid line}) corresponding to degree distributions with exponent $\beta=2.6$ ($\alpha=0.6$). 
In this case, the component size distributions can not feature a heavy tail, however, depending upon the scale parameter $s$ a transient asymptote with exponent $-1.6$ (\emph{dashed line}) emerges:
\emph{1}) $s=8.3\cdot 10^{-2},$ fast convergence to the exponential asymptote. 
 \emph{2}) $s=8.3\cdot 10^{-5},$ the distribution transiently follows what seems to be a heavy tail for $n<10^{6},$ whereas for larger $n$ the theory predicts no heavy tail. 
}
\label{fig:a06_C005_C0002}
\end{center}
\end{figure}

According to the definition \eqref{eq:Zol_mu}, the location parameter vanishes, $\mu(n)\equiv 0,$ when $\alpha < 1$. In this case, $x(n)=\frac{n}{\gamma(n)} \to 0 $ as $n\to \infty,$ and only \emph{one} asymptotic mode is possible for $w(n)$. 
Stable law $G^A(x,\alpha,1)$ is supported on $(0, \infty),$ and we make use of the series expansion around  $x\to0^+,$
\begin{multline}
\label{eq:lead_term2}
G^A(x,\alpha,1) =\\
 \frac{e^{-(1-\alpha) \left(\frac{\alpha}{x }\right)^{\frac{\alpha}{1-\alpha}}} \left(\frac{\alpha}{x }\right)^{\frac{1}{2} \left(1+\frac{1}{1-\alpha }\right)}}{\sqrt{2 \pi  \alpha(1-\alpha) }} 
[1+ O(x^{\frac{1-\alpha}{\alpha}})], \; x\to0^+,
\end{multline}
which when plugged in \eqref{eq:law} yields faster then algebraic decay of the component size distribution, see  Table~\ref{tab:asymptotes}, Case F.
Due to the parametrisation scheme for the stable law, the point $\alpha=1$ needs to be considered separately. In this case, $\frac{n-\mu(n)}{\gamma(n)}\to\infty$ when $n\to \infty,$ and we utilise the leading term of the series expansion,
\begin{equation}
\label{eq:lead_term3}
G^A(x,\alpha,1) =\frac{1}{\sqrt{2 \pi}}e^{\frac{x-1}{2} -e^{x-1}} [1+O(e^{1-x})],\; x\to\infty,
\end{equation}
which admits one sub-algebraic asymptotic mode for $w(n)$ as shown in Table~\ref{tab:asymptotes},~Case E.
This case is special in that the stable law $G^{A}(x,\alpha,1)$ is supported on $x\in(-\infty,\infty),$ but asymptotically, $\frac{1}{\gamma(n)}G^{A}\left(\frac{n-\mu(n)}{\gamma(n)},\alpha,1\right)$ always tends to $-\infty$ for large $n.$ At the same time, if for small $n$ the point $x(n)=\frac{n-\mu(n)}{\gamma(n)}$ stays on the positive half-axis where \eqref{eq:lead_term3} does not provide correct description for $G^{A}(x,\alpha,1),$   the convergence to the asymptote will be slow. In other words, there is an intermediate asymptote that the component size distribution can be approximated with, before it eventually switches to Eq.~\eqref{eq:lead_term3}. This switching point is given by such $n_0$ that $x(n)$ changes the sign from $'+'$ to $'-',$ i.e. when $n$ becomes greater then $n_0.$ By solving $x(n_0)=0,$ one obtains $n_0=e^{\frac{1}{s}}$, which means that in principle, the switching between the intermediate and the final asymptotes may be indefinitely  postponed if $s$ is small enough. The intermediate asymptote itself is deduced from the leading term of the stable law expansion at $\infty,$ that is $G^A(x,1,1)=\frac{2}{\pi} x^{-2}+O(x^{-3}).$ After the substitutions one obtains, 
$$w(n)\simeq\frac{\mu_1 s}{ (s \log n-1)^2}n^{-2},\; \frac{1}{s}>>0, \;n<e^{\frac{1}{s}},\; \alpha =1.$$
As illustrated in Fig.~\ref{fig:a06_C005_C0002}, similar considerations  are also valid for the case $0<\alpha<1,$ where
$$w(n)\simeq \mu_1 s   \frac{\Gamma(\alpha+1)}{\Gamma(\alpha)} n^{-\alpha-1} ,\; \frac{1}{s}>>0, \;0<\alpha <1.$$

When occurs, such switching has a practical importance when dealing with empirically observed component size data.
Indeed, it may happen that one observes only the intermediate asymptote and not the final one due a small number of samples at the tail of the component size distribution. For instance, the second curve in Fig.~\ref{fig:a06_C005_C0002} does feature an exponential decay at infinity, but if one limits the data points to $n < 10^6,$ the component size distribution will seem to be a  heavy-tailed one.

\begin{figure}[h]
\begin{center}
\includegraphics[width=0.45\textwidth]{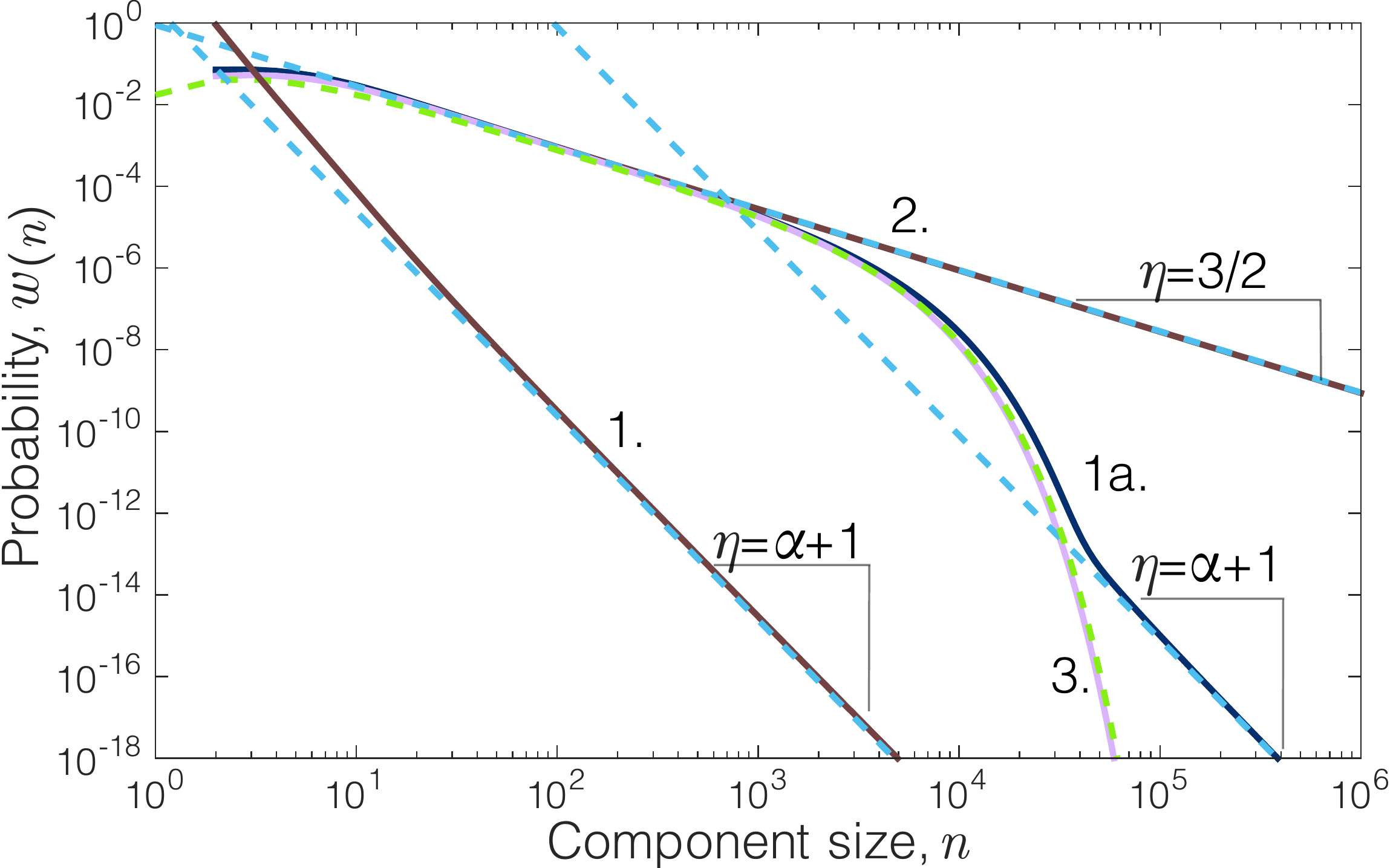}
\caption{Examples of component size distributions (\emph{solid lines}) that are associated with heavy-tailed degree distributions with $\beta=6$ ($\alpha=4$). The dashed lines represent the asymptotes in accordance with Case B in Table~\ref{tab:asymptotes}. 
Depending on the sign of $\theta,$ three asymptotic modes are distinguished: 
\emph{1)}  $s=1.93:$ $\theta=-0.8$, 
\emph{1a)} $s=9.42:$ $\theta=-0.027,$ 
\emph{2)}  $s=9.69:$ $\theta=1.6\cdot10^{-7}$, 
\emph{3)}  $s=10.05:$ $\theta=0.038.$ 
When $\theta$ is a small negative number (\emph{curve 1a}), $w(n)$ first decays as $n^{-3/2}$ but eventually switches to asymptote $n^{-\alpha-1}=n^{-5}.$  }
\label{fig:alpha4}
\end{center}
\end{figure}
Finally, we consider the case when the condition \eqref{eq:tail_cond} holds for $\beta>4:$ even though $u_1(k)$ has finite mean $M$ and variance $\sigma^2$ it also features a heavy  tail.
Again, as $n\to \infty,$ $x(n)=\frac{n-\mu(n)}{\sigma(n)}$ features the limiting values that are defined by the sign of $\theta$, see Eq.~\eqref{eq:point_asympt}.
One would expect that since $\sigma^2$ is finite, this case should be also well approximated with \eqref{eq:w_asm2}. This is indeed the case for $x(n)\approx 0.$ 
However, large deviations from zero $x(n)>>0$ do not follow Gaussian statistics \cite{ramsay2006,nacher2011}, and we approximate  $u_1^{*n}(k)$ with the Pareto stable law $ u_1^{*n}(k) \to \frac{1}{\sigma(n)} G^{P}(\frac{k-\mu(n)}{\sigma(n)},\alpha),\;n\to\infty.$ 
It turns out that $G^{P}(x,\alpha)$ behaves as the normal distribution for $x<C,$ where $C$ is a finite positive constant, but features a heavy tail with the same exponent as $u_1(k)$ when $x\to \infty,$ see Ref. \cite{ramsay2006}   
Thus, when $\theta \geq 0$ the component size distribution features asymptotic modes as in \eqref{eq:w_asm2},
while when $\theta < 0$ it features a heavy tail with exponent $-\alpha-1,$ see~Table~\ref{tab:asymptotes}, Case B.
Interestingly, when $\theta$ is a small negative number, $w(n)$ transiently follows one asymptote and then switches to the other as demonstrated in Figure~\ref{fig:alpha4}.
If there is a process that continuously changes the degree distribution so that $\theta$ progresses from being negative to positive, the exponent of the associated component size distribution will jump from the sub-critical branch, at $\theta<0$, to the critical one at $\theta=0.$ An example of such transition between two power-law modes is given in Figure~\ref{fig:alpha4}, where a component size distribution switches between power laws with exponents $\eta=3$ and $\eta=1.5$.

\section{Discussion and Conclusions}
 The broad generality of the results obtained in the previous section is achieved due to the fact that the configuration networks are locally tree-like and have vanishing probability of clustering in the thermodynamic limit, which allows one to benefit from the available in analytic combinatorics tools.
Eq. \eqref{eq:Lagrange1d}, that was analysed in the previous section, connects the degree distribution in a configuration network to the distribution of sizes for connected components.
The main conclusion one may draw from this equation is that the convolution power provides a smoothing effect. This means that all points of $u(k),\;k=1,\dots,\infty$ have a significant contribution to the definition of $w(n)$, but as $n$ increases, the system `forgets' the exact shape of the degree distribution and the component size distribution tends to the asymptote, that is defined by only a few parameters. The only information that is still preserved at the limit $n\to\infty$ is the first three moments of the degree distribution if such does not feature a heavy tail, see for example Fig.~\ref{fig:zigzag}. If $u(k)$ does feature a heavy tail then the information that characterises the tail becomes also important: that is the scale parameter $s$ and the exponent $\beta$. Depending upon the values of these parameters, many asymptotical modes exist. 

The expression for the asymptote is framed in terms of small deviation statistics for a sum of random variables and in some cases can be used as a good approximation for the component size distribution.
Table~\ref{tab:asymptotes} contains the analytical expressions for the asymptotes. Additionally, supporting code computing the component size distribution and the corresponding asymptotes is provided\cite{git:GECS}.
When using the asymptotical expressions to approximate $w(n)$, one should pay attention to two factors that follow from central limits: firstly $n$ should be large, secondly the approximation is best for $\theta$ close to zero. Finally, small deviations or a cutoff in a heavy-tailed degree distribution can trigger considerable and non-trivial changes in $w(n)$, for instance, the change of the asymptotical mode of the latter. 

\subsection{Degree distributions with a cutoff}
In practice, no empirical degree distribution is a heavy-tailed one. Most of the `real-world' degree distributions feature a cutoff, $u(k)=0,\; k>k_{\text{cut}},$ and therefore fail to be heavy-tailed in the strict sense of the definition \eqref{eq:tail_cond}. 
It turns out that if a cutoff is featured at large enough $k_{\text{cut}},$ the above-provided asymptotic analysis still has a relevant meaning. 
This situation can be compared to how we commonly attribute the fractal dimension to real-world geometric objects that fail to be  fractals on infinitesimal scales.

\begin{figure}[htbp]
\begin{center}
\includegraphics[width=0.4\textwidth]{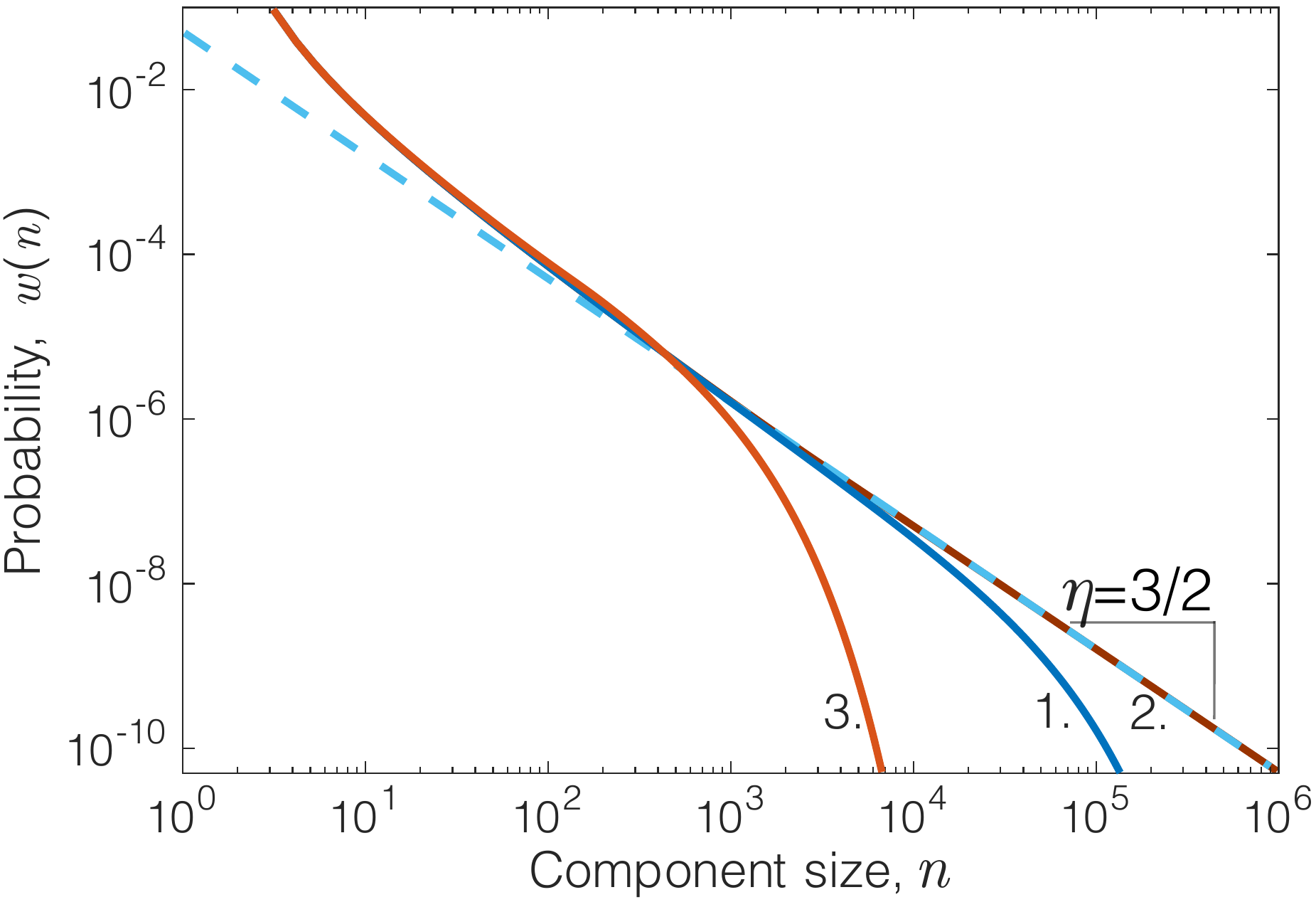}
\caption{The effect of a cutoff imposed on a heavy-tailed degree distribution with $\beta=3.3$ ($\alpha=1.3$) and $s=7.73.$ The solid curves correspond to component size distributions with: \emph{1)} no cutoff, $\theta>0$; \emph{2)} cutoff at $k=1000,$ $\theta=0$; \emph{3)} cutoff at $k=100,$ $\theta<0.$ The asymptote for \emph{1} is covered by Case D,   Table~\ref{tab:asymptotes}; while due to the cutoffs the asymptotes for \emph{2} and \emph{3} are covered by Case A.}
\label{fig:a13}
\end{center}
\end{figure}

Suppose one applies a cutoff at $k_{\text{cut}}$ to a degree distribution, $u(k),$  that features a heavy tail. Since $u(k)$ has a finite support, the asymptote of associated $w(n)$ is covered by Case A (Table~\ref{tab:asymptotes}), however, if $k_{cut}$ is large, $w(n)$ may also transiently follow the original asymptote.  Instead of an analytical investigation, we demonstrate the influence of the cutoff with numerical examples obtained by computing \eqref{eq:Lagrange1d}. This influence strongly depends on how the sign of $\theta$ is affected by the introduction of the cutoff. For example, if $\theta>0$ even after the cutoff, the cutoff will cause more nodes to appear in finite-size components, and thus the component size distribution will shift towards larger sizes. The opposite case is valid when $\theta\leq 0$ before (and after) the cutoff, then the cutoff causes the component size distribution to shift towards smaller sizes. The third option is when the cutoff changes the sign of $\theta$ form `$+$' to `$-$'. In this case, both shifts are possible. Fig.~\ref{fig:a13} shows how a component size distribution that corresponds to degree exponent $\beta=3.3$ is affected by a cutoff with various vales of $k_{\text{cut}}$.
\begin{figure}[htbp]
\begin{center}
\includegraphics[width=0.4\textwidth]{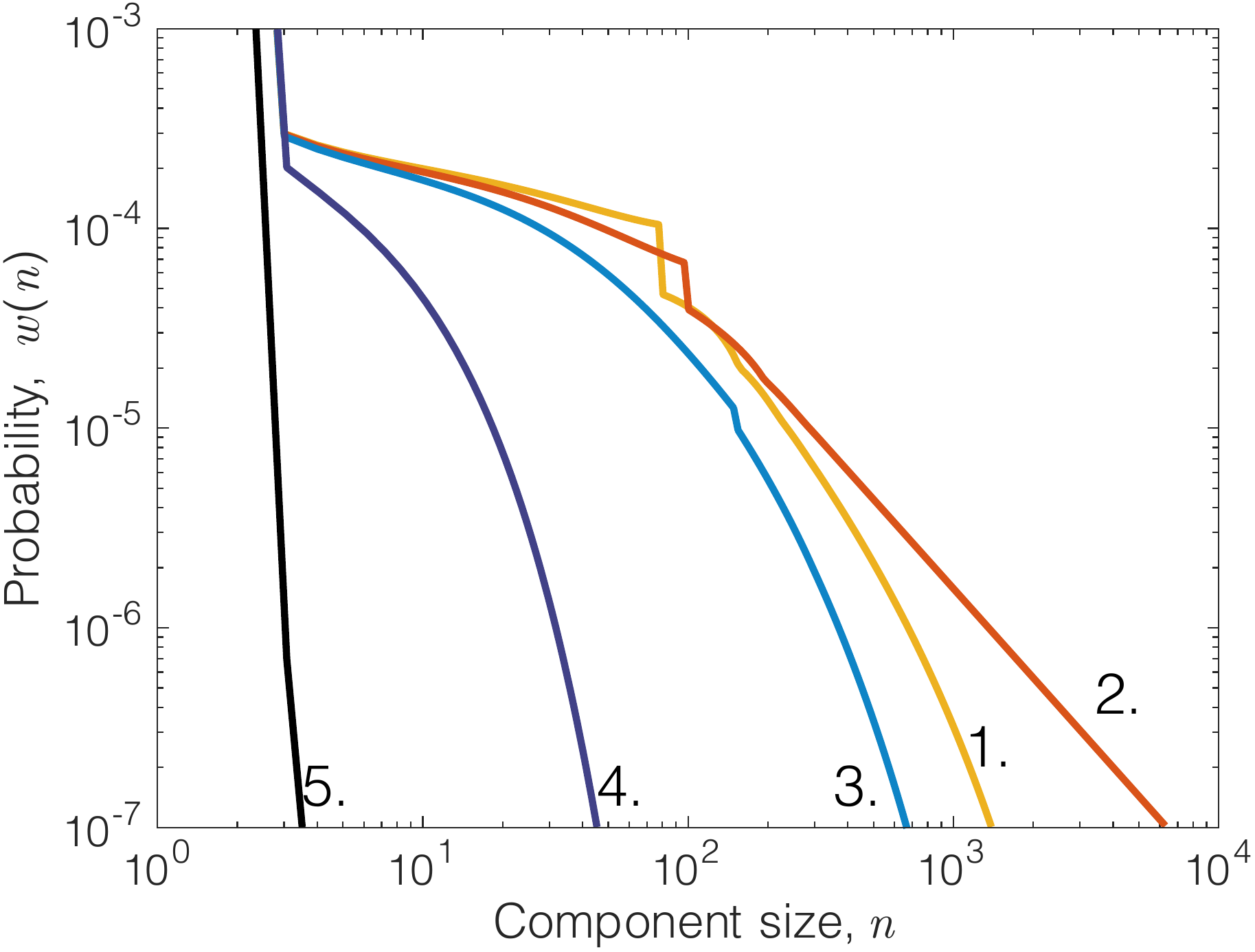}
\caption{The effect of a cutoff imposed on a heavy-tailed degree distribution with exponent $\beta=1$ ($\alpha=-1$) and $s=-2\cdot 10^{-4}$. The following values of the cutoff are considered: \emph{1}) $k_{\text{cut}}=80$, corresponds to $\theta<0$; \emph{2}) $k_{\text{cut}}=100,$ corresponds to $\theta \approx 0; \emph{3})$  $k_{\text{cut}}=150$ corresponds to $\theta>0,$ \emph{4}) $k_{\text{cut}}=10^3;$ \emph{5}) $k_{\text{cut}}=10^5.$
 }
\label{fig:a1}
\end{center}
\end{figure}
\subsection{Excess degree distribution with no mean value}
In principle, the excess degree distributions that do not have a mean value, i.e. $\beta<2$, do not fall within any of the above categories. However, if one introduces a cutoff, $u(k)$ will feature finite moments including, $\mu_3<\infty,$ hence this case should be treated according to Case A of Table~\ref{tab:asymptotes}.
Fig.~\ref{fig:a1} shows how cutoffs at $k=k_{\text{cut}}$ influence an instance of component size distribution with $\beta=1$. Unlike as in the previous example, in which $u(k)$ with no cutoff generates a valid $w(n)$, here the increase of $k_{\text{cut}}$ results in vanishing probability of finding a finite-size component at all: for any $n,$ $w(n)\to0$ when $k_{\text{cut}}\to\infty $. This illustrates the fact that finite-size components do not exist for $\beta\leq1,$ and the whole configuration network is connected almost surely. Non-existence of finite components for $\beta\leq 1$ also follows from the fact that in this case $\mu_1$ diverges and the point values of $w(n),$ as given below the the definition \eqref{eq:Lagrange1d*}, tend to zero.

Suppose the cutoff in the empirical, heavy-tailed degree distribution is due to the fact that the network sample has a finite size, $k_{cut}=N,\; N \neq \infty$, then one may approximate the expected number of edges in this sample as  
$$n_e =\frac{ N \mu_1}{2} = \frac{ N }{2} \sum_{k=1}^{N} k u(k) \simeq \frac{ N }{2} \sum_{k=1}^{N}  k^{-\beta+1}, \; N>>1,$$
so that
$$
n_e \simeq \begin{cases}
 N (1 - N^{2 - \beta } ), & \beta\neq 2,\\
N \log N, & \beta = 2.
\end{cases}
$$ 
Subsequently, three scenarios are possible here:\\
i)   \emph{sparse network,} $n_e=C N,\; C>0,\; \beta>2: $ the asymptotic modes are given in Table~\ref{tab:asymptotes};\\
ii)  \emph{semi-dense network,} either $n_e=C N \log N,\; \beta = 2 $ or $n_e=C N^{3-\beta},\; 1<\beta < 2:$ the mean value of excess distribution diverges; there are finite components but no power law in the distribution of component sizes;\\
iii) \emph{dense network, } $n_e=C N^{3-\beta},\; \beta \leq 1:$ the mean value of degree distribution $\mu_1\to \infty$, and finite components vanish as $N\to \infty$.

\begin{figure}
\begin{center}
\includegraphics[width=0.35\textwidth]{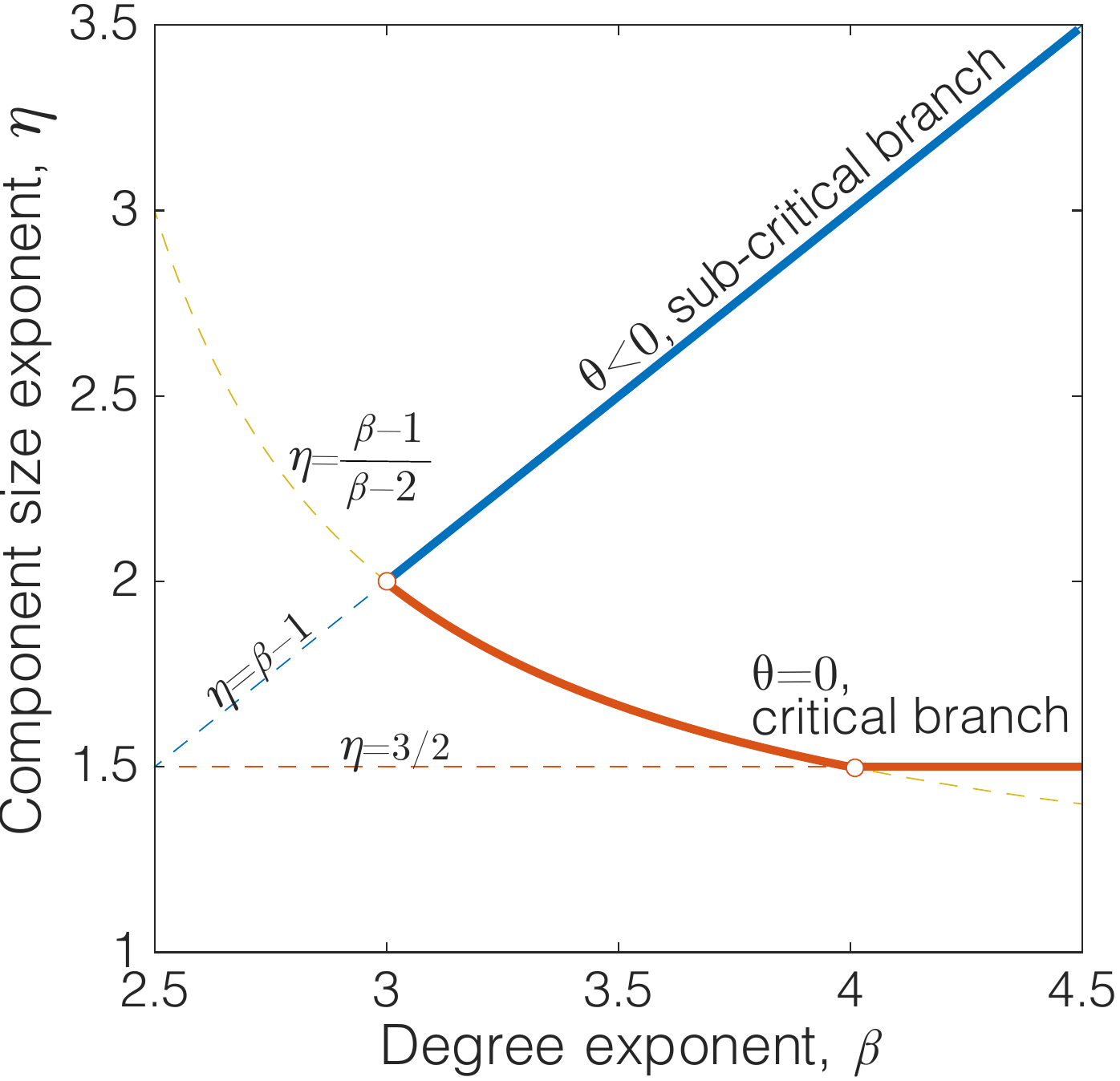}
\caption{The correspondence between the exponent $\beta$ in a heavy-tailed degree distribution, and the exponent $\eta$ in the associated heavy-tailed component size distribution. The critical brunch corresponds to $\theta=0,$ subcritical branch to $\theta<0.$ Positive $\theta$ is not associated with heavy-tailed component size distributions.}
\label{fig:diagram}
\end{center}
\end{figure}

\subsection{The role of the giant component}
All the cases presented in Table~\ref{tab:asymptotes} depend in some way on the value of $\theta$.
This is not a coincidence as the sign of $\theta$ is the indicator for the giant component existence. 
If the degree distribution features a heavy tail with exponent $\beta\geq3$, depending upon the value of $\theta,$ there are two possible heavy-tail exponents for the component size distribution: subcritical branch $\eta=\beta-1$ when $\theta<0$, and critical branch $\eta=\min\{\frac{3}{2},\frac{\beta-1}{\beta-2}\}$ when $\theta=0.$ 
This relation is illustrated in Fig.~\ref{fig:diagram}, where the component size distribution exponent $\eta$ is plotted versus the degree-distribution exponent $\beta$.
We can see that if the giant component exists, $\theta>0,$ then irrespectively of what is the degree distribution, the component size distribution always decays faster then the power law. Therefore it can be concluded that the giant component is not compatible with a heavy-tailed component size distribution. 
Any degree distribution with $\beta<3,$ leads to a giant component since $\theta$ can only be  positive in this case.
Furthermore, if $\beta\leq1$ the giant component is also the only component: with probability 1 the configuration network is fully connected.

\begin{table*}[htp]
\begin{center}
\begin{tabular}{|c| l|c|c|}
\hline
Finite moments of $u(k)\;\;$ 		 				 & $u(k)$, $k \to\infty$							& $\theta=\mu_2-2 \mu_1$ & Asymptote of $w(n)\;\;$	                          \\
\hline
\multirow{ 5}{*}{$\mu_3<\infty$}			 & \multirow{ 2}{*}{A.  $o(k^{-\beta}),\; \beta>4$ } 				 	& $\theta \neq 0$  					 & 	$C_1e^{-C_2 n} n^{-3/2}$ \\
																										
								 & 							 											& $\theta = 0$	 	 & 	$C_1 n^{-3/2}$								  \\ \cline{2-4}
						 		 & \multirow{ 3}{*}{B.  $O(k^{-\beta}),\; \beta>4$} 			& $\theta < 0$  	 & 	$C_3 n^{-\alpha-1}$					  \\
 								 & 																		& $\theta = 0$  	 &	$C_1 n^{-3/2}$			 					  \\ 
								 						 		 & 										& $\theta > 0$  	 &	$C_1 e^{-C_2 n}n^{-3/2}$ 					 \\ \hline	
 \multirow{ 6}{*}{
$\mu_3=\infty,$
$ \mu_2<\infty$
 } & \multirow{ 3}{*}{C. $O(k^{-\beta}),\; \beta=4$ }		& $\theta < 0$		 & $C_3 n^{-\alpha-1}$					  	        \\ 
&																								 		& $\theta = 0$		 &  $C_1'\frac{n^{-3/2} }{\sqrt{ \log n}}$                   						\\
 & 																										& $\theta > 0$ 		 & $C_1'\frac{n^{-3/2} }{\sqrt{ \log n}}e^{ -C_2'\frac{n}{\log n}}$	\\ 					  
	 \cline{2-4}
  &\multirow{ 3}{*}{D. $O( k^{-\beta}),\; 3<\beta<4$ } 																									& $\theta < 0$  	 &  $C_3 n^{-\alpha -1}$							\\
  & 																									& $\theta = 0$  	 &  $C_4 n^{-\frac{1}{\alpha }-1} $					\\ 	
 & 																										& $\theta > 0$  	 &  $C_5 e^{-C_6 n} n^{-3/2}$						\\ \cline{1-4}
 \multirow{ 2}{*}{$\mu_2=\infty$}   &\multirow{ 1}{*}{E. $O( k^{-\beta}),\; \beta=3$ }    	&  $\theta > 0$ 		 &  $C_7  e^{ -C_8 - C_9 n^{ \frac{2}{\pi} } }  n^{ \frac{1}{\pi} -2}$ 											 	\\ \cline{2-4}
  &\multirow{ 1}{*}{F. $O( k^{-\beta}),\; 2<\beta<3$ }    									&$\theta > 0$  	 &   $C_{10} e^{-C_{11} n} n^{-3/2}$						\\
 \hline
\end{tabular}
\begin{tabular}{|ll|}	
$C_1 = \frac{\mu_1^2}{ \sqrt{2 \pi (\mu_1 \mu_3-\mu_2^2) }}$, $\;\;C_1'=\frac{\mu_1}{\sqrt{2 \pi s}} ,$ & $C_7 = \frac{\sqrt{2}\mu_1}{\pi^{3/2} s}$,\\
$C_2 = \frac{(\mu_2-2 \mu_1)^2}{2( \mu_1 \mu_3- \mu_2^2)}$,  $\;\;C_2'=\frac{ (\mu_2-2\mu_1)^2}{2 s \mu_1^2} ,$ & $C_8 =\frac{1}{ \pi s}+\frac{1}{2}$,\\
$C_3 = \frac{s\mu_1^{\alpha+2}\Gamma(\alpha+1)}{(2\mu_1-\mu_2)^{\alpha+1}\Gamma(\alpha)}, $ & $C_9 = e^{-1 - \frac{2}{ \pi s}}$,\\
$C_4 = \mu_1\Gamma \left(1+\frac{1}{\alpha }\right) \sin \frac{\pi }{\alpha } \left(\frac{2 \Gamma (\alpha )\sin \frac{\pi  \alpha }{2} }{\pi ^{\alpha +1} s}\right)^{1/\alpha }\!,$ & $   C_{10}=\frac{\mu_1}{\sqrt{2-2\alpha }} \left(\frac{\sqrt{2}\Gamma (\alpha )\sin \frac{\pi  \alpha }{2} }{ \alpha  \pi ^{\alpha } s}\right)^{\frac{1}{2 \alpha -2}},$ \\
$    C_5 =\frac{\mu_1}{\sqrt{\alpha -1}} \left(\frac{2^{2-\alpha }\left(\frac{\mu_2}{\mu_1}-2\right)^{2-\alpha }  \Gamma (\alpha ) \sin \frac{\pi  \alpha }{2} }{\alpha  \pi ^{\alpha } s}\right)^{\frac{1}{2 \alpha -2}}\!\!\!\!,\;\;$& $   C_{11}=(1-\alpha ) \left(\frac{\sqrt{2}  \Gamma (\alpha )\sin \frac{\pi  \alpha }{2}}{\pi  \alpha ^{\alpha } s}\right)^{\frac{1}{\alpha -1}},\;\;\;\;$\\
$    C_6 = (1-\alpha ) \left(\frac{2 \left(\frac{\mu_2}{\mu_1}-2\right)^{\alpha } \Gamma (\alpha )\sin \frac{\pi  \alpha }{2} }{\alpha ^{\alpha }  \pi  s}\right)^{\frac{1}{\alpha -1}},$&$\alpha=\beta-2.$\\
\hline
\end{tabular}

\end{center}
\caption{Asymptotic behaviour of component sizes $w(n),$ in terms of degree distribution parameters: the first three moments $\mu_1,\mu_2,\mu_3$, scale parameter $s$, and exponent $\beta$. Supporting source code available in Ref.~\cite{git:GECS}.
 }
\label{tab:asymptotes}
\end{table*}%
 
\begin{acknowledgments}
This work is part of the research program Veni with project number 639.071.511, which is financed by the Netherlands Organisation for Scientific Research (NWO).
\end{acknowledgments}

\bibliographystyle{apsrev4-1}
\bibliography{literature}

	\end{document}